\documentclass[11 pt]{amsart}
\usepackage{amsmath, amsfonts, amsthm, color,latexsym}
\usepackage{graphicx}%
\providecommand{\U}[1]{\protect\rule{.1in}{.1in}}
\newtheorem{theorem}{Theorem}

\newtheorem{lemma}[theorem]{Lemma}

\newtheorem{problem}[theorem]{Problem}
\newtheorem{proposition}[theorem]{Proposition}
\newtheorem{remark}[theorem]{Remark}

\begin{document}
	
\title{Remarks on the Bohnenblust--Hille inequalities}
\author[Djair Paulino, Daniel Pellegrino and Joedson Santos]{Djair Paulino, Daniel Pellegrino and Joedson Santos}
\thanks{Djair Paulino is partially supported by Capes, Daniel Pellegrino is supported by CNPq Grant 307327/2017-5 and Grant 2019/0014 Paraíba State Research Foundation (FAPESQ) and Joedson Santos is supported by CNPq Grant 309466/2018-0 and Grant 2019/0014 Paraíba State Research Foundation (FAPESQ).\thinspace\hfill}
\thanks{\indent 2010 Mathematics Subject Classification: Primary 47A63; Secondary 47H60.}
\thanks{\indent Keywords: Bohnenblust--Hille inequalities; multilinear forms; sequence spaces.}
\date{}
\maketitle

\begin{abstract}
We revisit the Bohnenblust--Hille multilinear and polynomial inequalities and
prove some new properties. Our main result is a multilinear version of a recent
result on polynomials whose monomials have a uniformly bounded number of variables.
\end{abstract}


\section{Introduction}

Let $\mathbb{K}$ denote the real $\mathbb{R}$ or complex $\mathbb{C}$ scalar
fields. By $c_{0}$ we denote the Banach space of all scalar-valued sequences
$(a_{j})_{j=1}^{\infty}$ such that $\lim\limits_{j\longrightarrow\infty}%
a_{j}=0,$ endowed with the sup norm. For an $m$-linear form $T:c_{0}%
\times\cdots\times c_{0}\longrightarrow\mathbb{K}$ we denote, as usual,
\[
\Vert T\Vert:=\sup\{|T(x^{(1)},...,x^{(m)})|:\Vert(x^{(j)})\Vert
=1~~\mbox{for all}~~j=1,...,m\}.
\]
\textit{Littlewood's }$4/3$\textit{ inequality} \cite{LLL} (1930) asserts
that
\[
\left(  \sum_{j,k=1}^{\infty}\left\vert T(e_{j},e_{k})\right\vert ^{\frac
	{4}{3}}\right)  ^{\frac{3}{4}}\leq\sqrt{2}\Vert T\Vert,
\]
for all continuous bilinear forms $T\colon c_{0}\times c_{0}\rightarrow
\mathbb{K}$, and the exponent $4/3$ is sharp.
The constant $\sqrt{2}$ is sharp for real scalars.
Littlewood's $4/3$ inequality was perhaps one of the main cornerstones of the
beginning of several important inequalities on multilinear forms, such as an
inequality due to Bohnenblust and Hille (1931). The \textit{Bohnenblust--Hille
	inequality} \cite{bh} says that there is an optimal constant $B_{m}%
^{\mathbb{K}}\geq1$ such that
\[
\left(  \sum_{j_{1},\ldots,j_{m}=1}^{\infty}\left\vert T(e_{j_{1}}%
,\ldots,e_{j_{m}})\right\vert ^{\frac{2m}{m+1}}\right)  ^{\frac{m+1}{2m}}\leq
B_{m}^{\mathbb{K}}\left\Vert T\right\Vert ,
\]
for all continuous $m$-linear forms $T\colon c_{0}\times\cdots\times
c_{0}\rightarrow\mathbb{K}$. The optimal values of the constants
$B_{m}^{\mathbb{K}}$ are still unknown but when truncated to $\ell_{\infty}^{n}\times \cdots \times \ell_{\infty}^{n}$, the optimal values of the respective constants are computable within a finite number of elementary steps (see \cite{fv} and the references therein). The optimal values of these constants play a fundamental role in
applications in combinatorics and Quantum Information Theory (see \cite{Ara,
	MIT, montanaro} and the references therein).

Let $\alpha=(\alpha_{j})_{j=1}^{\infty}$ be a sequence in $\mathbb{N}%
\cup\{0\}$ and define $|\alpha|:=\sum\alpha_{j};$ we also denote $x^{\alpha
}:=\prod_{j}x_{j}^{\alpha_{j}}$. A mapping  $P:c_{0}%
\rightarrow\mathbb{C}$ is a continuous $m$-homogeneous polynomial if there exists a continuous $m$-linear form $\hat{P}: c_{0}\times \cdots \times c_{0} \longrightarrow \mathbb{C}$ such that for every $x\in c_{0}$ we have $P(x)= \hat{P}(x,...,x)$. For each sequence $\alpha$ satisfying  $|\alpha|=m,$ we denote by $c_{\alpha} (P)$ the coefficient of the monomial $x_{i_1}^{\alpha_{i_1}}\cdots x_{i_m}^{\alpha_{i_m}}$. Defining the norm of $\|P\|$ by $\|P\|:=\sup_{x\in B_{c_0}} |P(x)|,$ the Bohnenblust--Hille inequality for  $m$-homogeneous polynomial reads as follows: there is an optimal constant $C_{m}^{\mathbb{K}}\geq1$ such that
\[
\left(  \sum_{|\alpha|=m}|c_{\alpha}(P)|^{\frac{2m}{m+1}}\right)  ^{\frac
	{m+1}{2m}}\leq C_{m}^{\mathbb{K}}\Vert P\Vert
\]
for all continuous $m$-homogeneous polynomials $P:c_{0}\longrightarrow
\mathbb{K}$.

In \cite[Corollaries 3.2 and 3.3]{BAYART} it was proved that there exists
$C>0$ such that
\begin{equation}
B_{m}^{\mathbb{C}}\leq Cm^{\frac{1-\gamma}{2}}~~~\text{and}~~~~~~~B_{m}%
^{\mathbb{R}}\leq Cm^{\frac{2-log2-\gamma}{2}},\label{1112}%
\end{equation}
where $\gamma$ it is the Euler--Mascheroni constant. It is well-known that $C=1.3$ satisfies (\ref{1112}) (see \cite{Cavalcante}); numerically,
$\frac{1-\gamma}{2}\simeq0.211392$ and $\frac{2-\log2-\gamma}{2}\simeq
0.36482$. Recent advances suggest that the above estimates are far from being
sharp; it is even conjectured that for real scalars the optimal constants in
(\ref{1112}) are surprisingly $2^{1-1/m}$ and in the complex case it is
sometimes conjectured that the optimal constants may be the trivial ones (see
\cite{fv} and the references therein).

Still in \cite[Corollaries 5.3 and 5.4]{BAYART}, it has been proven that, for
any $\varepsilon>0$, there exists $\kappa>0$ such that, for any $m\geq1$,
\begin{equation}
C_{m}^{\mathbb{C}}\leq\kappa\left(  1+\varepsilon\right)  ^{m}.\label{2qqq}%
\end{equation}
Observing the upper bounds for the polynomial and multilinear
Bohnenblust--Hille inequalities, it seems natural to speculate that
$B_{m}^{\mathbb{K}}\leq C_{m}^{\mathbb{K}}$ for all $m.$ However, the best
known lower bounds for $C_{m}^{\mathbb{C}}$ are almost the trivial ones (see \cite{Nunez}) and from this perspective the inequality $B_{m}^{\mathbb{K}}\leq
C_{m}^{\mathbb{K}}$ appears less evident. As a matter of fact, we have found no information in this direction in the literature, and for this reason we now present a proof that the multilinear constant is dominated by the polynomial constant.

\begin{proposition}
	\label{compconst} The optimal constant of the multilinear Bohnenblust--Hille inequality is smaller than optimal constant of the polynomial Bohnenblust--Hille inequality, that is,
	\[
	B_{m}^{\mathbb{K}}\leq C_{m}^{\mathbb{K}}.
	\]
	
\end{proposition}

\textit{Proof.} Let $T\colon c_{0}\times\cdots\times c_{0}\rightarrow
\mathbb{K}$ be a continuous $m$-linear form. Let $\mathbb{N}=\mathbb{N}%
_{1}\cup\cdots\cup\mathbb{N}_{m}$ be a pairwise disjoint union, with
$card(\mathbb{N})=card(\mathbb{N}_{r}),$ for all $r=1,...,m$. Consider
bijections $\sigma_{i}:\mathbb{N}\longrightarrow\mathbb{N}_{i}$, $i=1,...,n$,
and define $T_{1}\colon c_{0}\times\cdots\times c_{0}\rightarrow\mathbb{K}$
by
\[
T_{1}((x_{j}^{1})_{j=1}^{\infty},...,(x_{j}^{m})_{j=1}^{\infty})=T((x_{\sigma
	_{1}(j)}^{1})_{j=1}^{\infty},...,(x_{\sigma_{m}(j)}^{m})_{j=1}^{\infty}).
\]
Note that $T_{1}$ has the same coefficients of $T$ and, therefore, $\Vert
T_{1}\Vert=\Vert T\Vert$. Define $P:c_{0}\longrightarrow\mathbb{K}$ by
$P(x)=T_{1}(x,...,x)$. So $\Vert P\Vert\leq\Vert T_{1}\Vert,$ and it follows
that
\begin{align*}
\left(  \sum_{j_{1},...,j_{m}=1}^{\infty}|T(e_{j_{1}},...,e_{j_{m}}%
)|^{\frac{2m}{m+1}}\right)  ^{\frac{m+1}{2m}}  &  =\left(  \sum_{j_{1}
	,...,j_{m}=1}^{\infty}|T_{1}(e_{j_{1}},...,e_{j_{m}})|^{\frac{2m}{m+1}%
}\right)  ^{\frac{m+1}{2m}}\\
&  =\left(  \sum_{|\alpha|=m}|c_{\alpha}(P)|^{\frac{2m}{m+1}}\right)
^{\frac{m+1}{2m}}\\
& \leq C_{m}^{\mathbb{K}}\Vert P\Vert\leq C_{m}^{\mathbb{K}}\Vert T\Vert.
\end{align*}
Therefore,
\[
B_{m}^{\mathbb{K}}\leq C_{m}^{\mathbb{K}}.
\]
\hfill$\Box$

\begin{remark}
	The same argument used in the proof of Proposition \ref{compconst} can be used
	to show that a similar phenomenon occurs with the Hardy--Littlewood
	inequalities. 
\end{remark}

In addition to the classical Bohnenblust--Hille inequalities for polynomials
and multilinear forms, some authors have studied these inequalities for some
restricted classes . For example, the multilinear Bohnenblust--Hille inequality
for unimodular coefficients (the coefficients are $\pm1$) was studied in
\cite{Teixeira}. A fundamental result on unimodular multilinear forms is the following (see \cite[Lemma 6.1]{ABPS}):

\begin{theorem}
	[Kahane--Salem--Zygmund inequality]\label{sss} Let $m,n\geq1.$ There is a constant $K_{m}>0$, depending only on $m$, and an $m$-linear form
	$T_{m,n}\colon c_{0}\times\cdots\times c_{0}\rightarrow\mathbb{K}$ of the
	form
	\[
	T_{m,n}(z^{(1)},...,z^{(m)})=\displaystyle\sum_{i_{1},...,i_{m}=1}^{n}\pm
	z_{i_{1}}^{(1)}\cdots z_{i_{m}}^{(m)}
	\]
	such that
	\[
	\Vert T_{m,n}\Vert\leq K_{m}n^{\frac{m+1}{2}}.
	\]
	
\end{theorem}

Since the form given by the Kahane--Salem--Zygmund (KSZ) inequality is
unimodular, it can be easily proved that the optimal exponent of the
Bohnenblust--Hille inequality
\begin{equation}
\left(  \sum_{i_{1},...,i_{m}=1}^{\infty}|T(e_{i_{1}},...,e_{i_{m}}%
)|^{\frac{2m}{m+1}}\right)  ^{\frac{m+1}{2m}}\leq B_{m}^{\mathbb{K}}\Vert
T\Vert\label{eqprincbh}%
\end{equation}
restricted to unimodular $m$-linear forms $T\colon c_{0}\times\cdots\times
c_{0}\rightarrow\mathbb{K}$ is also $\frac{2m}{m+1}$.

Polynomial versions of the KSZ inequality are natural consequences of the
respective multilinear version; however the coefficients of the polynomial KSZ inequalities are not unimodular
anymore (see \cite[Theorem 4]{boas}). However, according to \cite{Tomaz}, we have a polynomial version of Theorem \ref{sss} for polynomials with coefficients in $\{0,1,-1\}$. As an immediate consequence, that we were not able to find in the literature, we conclude that the exponent $\frac{2m}{m+1}$ of
Bohnenblust--Hille inequality
\[
\left(  \sum_{|\alpha|=m}|c_{\alpha}(P)|^{\frac{2m}{m+1}}\right)  ^{\frac
	{m+1}{2m}}\leq C_{m}^{\mathbb{K}}\Vert P\Vert
\]
restricted to $m$-homogeneous polynomials $P:c_{0}%
\longrightarrow\mathbb{K}$ with coefficients in  $\{0,1,-1\}$ is optimal.

Other constrained Bohnenblust--Hille inequalities were studied in \cite{Defant}
and \cite{Mariana}, for polynomials whose monomials have a uniformly bounded
number of variables.

Let us choose $card(A)$ to represent the cardinality of the set $A$. For positive
integers $m$ and $M\leq m$ we define
\[
\omega(\alpha)=card\{j:\alpha_{j}\neq0\}
\]
and
\[
\Lambda_{M,m}=\{\alpha:|\alpha|=m,~\omega(\alpha)\leq M\}.
\]

In \cite{Mariana} it was proved that for all integers $m$ and $M$, with $M\leq
m$, there exists a universal constant $k_{M}\geq1$ such that
\begin{equation}
\left(  \sum_{\alpha\in\Lambda_{M,m}}|c_{\alpha}(P)|^{\frac{2m}{m+1}}\right)
^{\frac{m+1}{2m}}\leq k_{M}\Vert P\Vert, \label{inequalitykM}%
\end{equation}
for all continuous $m$-homogeneous polynomials $P:c_{0}\longrightarrow
\mathbb{C}$, that is, it is possible to obtain a constant independently of the
value of $m$.

The main result of the present paper is a multilinear version of (\ref{inequalitykM}). It is worth observing that while the polynomial case is restricted to complex scalars, the multilinear result is of interest for the real and complex frameworks.

\section{Bohnenblust--Hille inequality with restrictions on the monomials}

Observing the proof of (\ref{inequalitykM}) in \cite{Mariana} we note that the
universal constant $k_{M}$ is obtained as a consequence of a stronger
inequality with exponent $\frac{2M}{M+1}$ instead of $\frac{2m}{m+1}.$ A
natural question is whether the exponent $\frac{2M}{M+1}$ is sharp. The next
result gives a lower bound for the optimal exponent.

\begin{proposition}
	\label{expotimo} Let $m$ and $M$ be positive integers, let $M\leq m$ and let
	$r\geq0$. If there is a constant $\delta_{M}^{\mathbb{C}}$ (not depending on
	$m$) such that
	\begin{equation}
	\left(  \sum_{\alpha\in\Lambda_{M,m}}|c_{\alpha}(P)|^{r}\right)  ^{\frac{1}
		{r}}\leq\delta_{M}^{\mathbb{C}}\Vert P\Vert, \label{desigualdadeexpo}%
	\end{equation}
	for all continuous $m$-homogeneous polynomials $P:c_{0}\longrightarrow
	\mathbb{C}$, then $r\geq\frac{2(M-1)}{M}$.
\end{proposition}

\textit{Proof.} Let $P:c_{0}\longrightarrow\mathbb{C}$ be a continuous
$(M-1)$-homogeneous polynomial. Define the $m$-homogeneous polynomial
$\overline{P}:c_{0}\longrightarrow\mathbb{C}$ by $\overline{P}(x)=x_{1}%
^{m-M+1}\cdot P(x).$ Note that each monomial of $\overline{P}$ has at most $M$
distinct indexes and $\Vert\overline{P}\Vert\leq\Vert P\Vert$. Therefore,
\[
\left(  \sum_{|\alpha|=M-1}|c_{\alpha}(P)|^{r}\right)  ^{\frac{1}{r}}=\left(
\sum_{\alpha\in\Lambda_{M,m}}|c_{\alpha}(\overline{P})|^{r}\right)  ^{\frac
	{1}{r}}\leq\delta_{M}^{\mathbb{K}}\Vert\overline{P}\Vert\leq\delta
_{M}^{\mathbb{K}}\Vert P\Vert.
\]
Since this inequality holds for any continuous $(M-1)$-homogeneous polynomial
$P:c_{0}\longrightarrow\mathbb{K}$, it follows from the optimality of the
exponent $\frac{2(M-1)}{M}$ of the polynomial Bohnenblust--Hille inequality
(for $(M-1)$-homogeneous polynomials) that $r\geq\frac{2(M-1)}{M}.$
\hfill$\Box$

The next result is a kind of multilinear version of the main result of \cite{Mariana}; as we already mentioned in the introduction, differently of the polynomial case, our result has also interest in the real case:

\begin{theorem}
	\label{Teoprincipal} Let $m,M$ be positive integers, with $M\leq m$. Then,
	there is a universal constant $\eta_{M}\geq1$ such that
	\begin{equation}
	\left(  \sum_{card(\{i_{1},...,i_{m}\})\leq M}|T(e_{i_{1}},...,e_{i_{m}}
	)|^{\frac{2m}{m+1}}\right)  ^{\frac{m+1}{2m}}\leq\eta_{M}\Vert T\Vert
	\label{eqprinc}%
	\end{equation}
	for all continuous $m$-linear forms $T\colon c_{0}\times\cdots\times
	c_{0}\rightarrow\mathbb{K}$. Moreover, if $\overline{C}_{m,M}^{\mathbb{R}}$
	is  the optimal constant for real scalars, then
	\begin{align*}
	2  &  \leq\liminf\limits_{m\rightarrow\infty}\overline{C}_{m,M}^{\mathbb{R}
	}\leq\limsup_{m\rightarrow\infty}\overline{C}_{m,M}^{\mathbb{R}}\leq
	M^{\frac{M+1}{2}}\text{, if }M\geq3,\\
	\sqrt{2}  &  \leq\liminf\limits_{m\rightarrow\infty}\overline{C}_{m,M}
	^{\mathbb{R}}\leq\limsup_{m\rightarrow\infty}\overline{C}_{m,M}^{\mathbb{R}
	}\leq\sqrt{8}\text{, if }M=2.
	\end{align*}
	
\end{theorem}

\begin{remark}
	It is interesting to observe that the above theorem is closely connected to
	the main result of \cite{AANPR} but, while in \cite{AANPR} the constants are
	contractive when $M$ is fixed and $m\rightarrow\infty$, in our result, for real scalars, the constants do not converge to $1.$
\end{remark}

Let us begin by establishing some notation. Let $n$ be a positive integer and
from now on $e_{i}^{n}$ denotes the $n$-tuple $\left(  e_{i},{...}%
,e_{i}\right)  $. Furthermore, if $n_{1},\ldots,n_{k}\geq1$ are such that
$n_{1}+\cdots+n_{k}=m$, then $\left(  e_{i_{1}}^{n_{1}},\ldots,e_{i_{k}%
}^{n_{k}}\right)  $ represents the $m$-tuple:
\[
(e_{i_{1}},\overset{\text{{\tiny $n_{1}$\thinspace times}}}{\ldots},e_{i_{1}%
},\ldots,e_{i_{k}},\overset{\text{{\tiny $n_{k}$\thinspace times}}}{\ldots
},e_{i_{k}}).
\]

The following theorem, recently proved in \cite{AANPR}, plays a fundamental
role in the proof of our main result.

\begin{theorem}
	[Bohnenblust--Hille inequality by blocks]\label{Lemabloco} Let $m\geq M\geq
	1$,  and let $n_{1},...,n_{M}\geq1$ be such that $n_{1}+\cdots+n_{M}=m$. Then,
	for  every continuous $m$-linear form $T:c_{0}\times\cdots\times c_{0}
	\longrightarrow\mathbb{K}$,
	\[
	\left(  \sum_{i_{1},...,i_{M}=1}^{\infty}|T(e_{i_{1}}^{n_{1}},...,e_{i_{M}
	}^{n_{M}})|^{\frac{2M}{M+1}}\right)  ^{\frac{M+1}{2M}}\leq B_{M}^{\mathbb{K}
	}\Vert T\Vert.
	\]
	Moreover, the exponent is optimal.
\end{theorem}

We also need the following consequence of the Khinchin inequality:

\begin{lemma}
	(see \cite[p. 14]{Teixeira}) \label{Lemainterpolar} For all continuous
	$m$-linear  forms $T:c_{0}\times\cdots\times c_{0}\longrightarrow\mathbb{K},$
	we have
	\[
	\left(  \sum_{i_{1},...,i_{m}=1}^{\infty}|T(e_{i_{1}},...,e_{i_{m}}
	)|^{2}\right)  ^{\frac{1}{2}}\leq\Vert T\Vert.
	\]
	
\end{lemma}

To prove the first part of Theorem \ref{Teoprincipal} we shall use a
combinatorial argument together with Theorem \ref{Lemabloco} and Lemma
\ref{Lemainterpolar} (that is valid for $\mathbb{K}=\mathbb{R}~~$%
or$~~\mathbb{C}$), combined with a kind of interpolation argument. To analyze
the constant when $m\rightarrow\infty$ we will use specific $m$-linear forms,
which have already been used in other works, as in \cite{Teixeira}.

Let $T:c_{0}\times\cdots\times c_{0}\longrightarrow\mathbb{K}$ be a continuous $m$-linear form and consider the set
\[
A_{m,M}:=\{(j_{i_{1}},...,j_{i_{m}}):i_{1},...,i_{m}\in\{1,...,M\}\}.
\]
Note that $card(A_{m,M})=M^{m}$. Now, for each choice $i_{1},...,i_{m}%
\in\{1,...,M\}$, consider the sum
\[
\sum_{j_{i_{1}},...,j_{i_{m}}=1}^{\infty}|T(e_{j_{i_{1}}},...,e_{j_{i_{m}}%
})|^{\frac{2M}{M+1}}.
\]
Note that
\[
\sum_{card(\{i_{1},...,i_{m}\})\leq M}\hspace{-0.1cm}|T(e_{i_{1}},...,e_{i_{m}})|^{\frac{2M}
	{M+1}}\leq\hspace{-0.1cm}\sum_{i_{1},...,i_{m}=1}^{M}\hspace{-0.1cm}\left(  \sum_{j_{i_{1}},...,j_{i_{m}}
	=1}^{\infty}\hspace{-0.1cm}|T(e_{j_{i_{1}}},...,e_{j_{i_{m}}})|^{\frac{2M}{M+1}}\hspace{-0.1cm}\right)  .
\]
By Theorem \ref{Lemabloco} (noting
that it is also valid if the same indices are not necessarily together), we
have
\begin{align*}
\sum_{card(\{i_{1},...,i_{m}\})\leq M}|T(e_{i_{1}},...,e_{i_{m}})|^{\frac{2M}
	{M+1}}  &  \leq\sum_{i_{1},...,i_{m}=1}^{M}(B_{M}^{\mathbb{K}})^{\frac
	{2M}{M+1}}\Vert T\Vert^{\frac{2M}{M+1}}\\
&  =M^{m}(B_{M}^{\mathbb{K}})^{\frac{2M}{M+1}}\Vert T\Vert^{\frac{2M}{M+1}}.
\end{align*}

Therefore,
\[
\left(  \sum_{card(\{i_{1},...,i_{m}\})\leq M}|T(e_{i_{1}},...,e_{i_{m}}%
)|^{\frac{2M}{M+1}}\right)  ^{\frac{M+1}{2M}}\leq M^{\frac{m(M+1)}{2M}}%
B_{M}^{\mathbb{K}}\Vert T\Vert.
\]
By Lemma \ref{Lemainterpolar}, we have
\[
\left(  \sum_{card(\{i_{1},...,i_{m}\})\leq M}|T(e_{i_{1}},...,e_{i_{m}}%
)|^{2}\right)  ^{\frac{1}{2}}\leq\left(  \sum_{i_{1},...,i_{m}=1}^{\infty
}|T(e_{i_{1}},...,e_{i_{m}})|^{2}\right)  ^{\frac{1}{2}}\leq\Vert T\Vert.
\]
Therefore, using interpolation (or the H\"older inequality for mixed sums)
with $\theta=M/m$, we have
\begin{align*}
\left(  \sum_{card(\{i_{1},...,i_{m}\})\leq M}|T(e_{i_{1}},...,e_{i_{m}}%
)|^{\frac{2m}{m+1}}\right)  ^{\frac{m+1}{2m}}  &  \leq\left(  M^{\frac
	{m(M+1)}{2M}}B_{M}^{\mathbb{K}}\Vert T\Vert\right)  ^{\theta}\cdot(\Vert
T\Vert)^{1-\theta}\\
&  =M^{\frac{M+1}{2}}(B_{M}^{\mathbb{K}})^{M/m}\Vert T\Vert.
\end{align*}
By (\ref{1112}), we have
\begin{align*}
\left(  \sum_{card(\{i_{1},...,i_{m}\})\leq M}|T(e_{i_{1}},...,e_{i_{m}}%
)|^{\frac{2m}{m+1}}\right)  ^{\frac{m+1}{2m}}  &  \leq M^{\frac{M+1}{2}%
}(1.3M^{0.365})^{M/m}\Vert T\Vert\\
&  =\left(  1.3\right)  ^{M/m}M^{\frac{0.365M}{m}+\frac{M+1}{2}}\Vert T\Vert.
\end{align*}

Let $\overline{C}_{m,M}^{\mathbb{R}}$ be the optimal constant of
inequality (\ref{eqprinc}) for real scalars. Consider, as well as in
\cite[page 15]{Teixeira}, the bilinear form $S_{2}:c_{0}\times c_{0}%
\longrightarrow\mathbb{R}$,
\[
S_{2}(x,y)=x_{1}y_{1}+x_{1}y_{2}+x_{2}y_{1}-x_{2}y_{2}.
\]

For $m=3$, consider $S_{3}: c_{0}\times c_{0} \times c_{0} \longrightarrow
\mathbb{R} $,
\begin{align*}
S_{3}(x,y,z)  &  =(z_{1}+z_{2})(x_{1}y_{1}+x_{1}y_{2}+x_{2}y_{1}-x_{2}y_{2})\\
&  ~~+(z_{1}-z_{2})(x_{3}y_{1}+x_{3}y_{2}+x_{4}y_{1}-x_{4}y_{2}).
\end{align*}
Note that, $S_{3}$ has at most three distinct indexes in each monomial and that $\|S_{3}\|=4.$

For $m=4$, define the form $S_{4}: c_{0}\times c_{0} \times c_{0} \times
c_{0}\longrightarrow\mathbb{R}$
\begin{align*}
S_{4}(x,y,z,w)  &  =(w_{1}+w_{2})((z_{1}+z_{2})(x_{1}y_{1}+x_{1}y_{2}%
+x_{2}y_{1}-x_{2}y_{2})\\
&  ~~+(z_{1}-z_{2})(x_{3}y_{1}+x_{3}y_{2}+x_{4}y_{1}-x_{4}y_{2}))\\
&  ~~+(w_{1}-w_{2})((z_{1}+z_{2})(x_{5}y_{1}+x_{5}y_{2}+x_{6}y_{1}-x_{6}%
y_{2})\\
&  ~~+(z_{1}-z_{2})(x_{7}y_{1}+x_{7}y_{2}+x_{8}y_{1}-x_{8}y_{2})).
\end{align*}
Note that, as previously, $S_{4}$ has at most three distinct indexes in each
monomial and that $\|S_{4}\|=8$.

The construction can be carried out by induction for all $m$-linear forms $S_m$ (see \cite[page 15]{Teixeira}), where in $S_{m}$ each monomial has at most three distinct indexes and, for all $m$, we have $$\Vert
S_{m}\Vert=2^{m-1}.$$

Note that
\[
\left(  \sum_{i_{1},...,i_{m}=1}^{\infty}|S_{m}(e_{i_{1}},...,e_{i_{m}%
})|^{\frac{2m}{m+1}}\right)  ^{\frac{m+1}{2m}}=(2^{2(m-1)})^{\frac{m+1}{2m}%
}=2^{\frac{(m-1)(m+1)}{m}}.
\]
Note that in the above argument we always have $M=3$. Therefore, considering
this form $S_{m}$ in (\ref{eqprinc}), we have

\[
2^{\frac{(m-1)(m+1)}{m}}\leq\eta_{M}2^{m-1},
\]
that is,
\[
\eta_{M} \geq 2^{\frac{(m-1)}{m}}.
\]
Thus
\[
\overline{C}_{m,M}^{\mathbb{R}}\geq2^{\frac{(m-1)}{m}}\text{, for all }%
M\geq3,
\]
and thus%
\[
\liminf\limits_{m\rightarrow\infty}\overline{C}_{m,M}^{\mathbb{R}}\geq2\text{,
	if }M\geq3.
\]

If $M=2$, we consider the $m$-linear forms (with $m$ even) given by%
\[
R_{m}(x^{(1)},...,x^{(m)})=%
{\textstyle\prod\limits_{i=1}^{m/2}}
\left(  x_{1}^{(2i-1)}x_{1}^{(2i)}+x_{1}^{(2i-1)}x_{2}^{(2i)}+x_{2}%
^{(2i-1)}x_{1}^{(2i)}-x_{2}^{(2i-1)}x_{2}^{(2i)}\right).
\]

It can be easily proved that%
\[
\left\Vert R_{m}\right\Vert =2^{m/2}%
\]
and, since $R_{m}$ has precisely $4^{m/2}$ monomials, replacing $R_m$ in (\ref{eqprinc}),  we conclude that%
\[
(4^{m/2})^{\frac{m+1}{2m}}  \leq \eta_{M} 2^{m/2}.
\]
So, $\eta_{M} \geq 2^{\frac{1}{2}}$ and
\begin{equation}\label{eqm21}
\liminf\limits_{k\rightarrow\infty}\overline{C}_{2k,2}^{\mathbb{R}}\geq\sqrt
{2}.
\end{equation}

If $m$ is odd, we consider the $m$-linear forms given by%

\begin{eqnarray*}
	&~&A_{m}(x^{(1)},..., x^{(m)})=\\
	&~&=%
	{\textstyle\prod\limits_{i=1}^{(m-1)/2}}
	x_1^{(m)}\cdot \left[  x_{1}^{(2i-1)}x_{1}^{(2i)}+x_{1}^{(2i-1)}x_{2}^{(2i)}+x_{2}%
	^{(2i-1)}x_{1}^{(2i)}-x_{2}^{(2i-1)}x_{2}^{(2i)}\right]  .
\end{eqnarray*}

It can be easily proved that%
\[
\left\Vert A_{m}\right\Vert =2^{(m-1)/2}%
\]
and, since $A_{m}$ has precisely $4^{(m-1)/2}$ monomials, replacing $A_m$ in (\ref{eqprinc}), we conclude that%
\[
(4^{(m-1)/2})^{\frac{m+1}{2m}}  \leq \eta_{M} 2^{(m-1)/2}.
\]
So, $\eta_{M} \geq 2^{\frac{m-1}{2m}}$ and
\begin{equation}\label{eqm22}
\liminf\limits_{k\rightarrow\infty}\overline{C}_{2k+1,2}^{\mathbb{R}}\geq\sqrt
{2}.
\end{equation}
Therefore, by (\ref{eqm21}) and (\ref{eqm22}), we have
\[
\liminf\limits_{m\rightarrow\infty}\overline{C}_{m,2}^{\mathbb{R}}\geq\sqrt
{2}.
\]
\hfill$\Box$

It is natural to wonder to what extent our main result can be improved. We
present below a couple of open problems arisen from our main result:

\bigskip

\begin{problem}
	What about the estimates of $\liminf\limits_{m\rightarrow\infty}\overline
	{C}_{m,M}^{\mathbb{K}}$ and $\limsup\limits_{m\rightarrow\infty}\overline{C}
	_{m,M}^{\mathbb{K}}$ for $\mathbb{K}=\mathbb{C}$?
\end{problem}

\begin{problem}
	\label{5}What are the optimal values of the real case $\overline{C}
	_{m,M}^{\mathbb{R}}$?
\end{problem}

Problem \ref{5} may be too complicated. Somewhat intermediate problems are:

\begin{problem}
	What about the optimal values of $\liminf\limits_{m\rightarrow\infty}
	\overline{C}_{m,M}^{\mathbb{R}}$ and $\limsup\limits_{m\rightarrow\infty}
	\overline{C}_{m,M}^{\mathbb{R}}$?
\end{problem}

\begin{problem}
	Is it possible to improve (asymptotically speaking) the upper estimate
	$M^{\frac{M+1}{2}}$?
\end{problem}

\textbf{Acknowledgement}. The authors would like to thank the referees for the detailed and insightful suggestions that helped fix some issues with the original version and make the text clearer for the readers.

\bigskip

\noindent Departamento de Matem\'atica\newline Universidade Federal da
Para\'iba\newline58.051-900 -- Jo\~ao Pessoa -- Brazil\newline e-mail:
djair.paulino@ifrn.edu.br and djairpsc@hotmail.com\newline 
e-mail: pellegrino@pq.cnpq.br and and dmpellegrino@gmail.com\newline
e-mail: joedson@mat.ufpb.br and joedsonmat@gmail.com

\frenchspacing

\begin{thebibliography}{99}                                       
	\bibitem {AANPR}N. Albuquerque, G. Ara\'{u}jo, W. Cavalcante, T. Nogueira, D.
	N\'{u}\~{n}ez--Alarc\'{o}n, D. Pellegrino, P. Rueda, On summability of
	multilinear operators and applications. \emph{Ann. Funct. Anal.} \textbf{9}
	(2018), no. 4, 574--590.
	
	\bibitem {ABPS}N. Albuquerque, F. Bayart, D. Pellegrino, J. B.
	Seoane--Sep\'ulveda, Sharp generalizations of the multilinear
	Bohnenblust--Hille inequality. \emph{J. Funct. Anal.} \textbf{266} (2014), 3726--3740.
	
	
	\bibitem {Ara}G. Ara\'{u}jo, D. Pellegrino, A Gale--Berlekamp
	permutation-switching problem in higher dimensions. \emph{European J. Combin.}
	\textbf{77} (2019), 17--30.
	
	\bibitem {MIT}S. Arunachalam, S. Chakraborty, M. Koucky, N. Saurabh, R. de
	Wolf, Improved bounds on Fourier entropy and Min-entropy. \emph{Electron.
		Colloquium Comput. Complexity} (ECCC) \textbf{25} (2018), 167.
	
	\bibitem {BAYART}F. Bayart, D. Pellegrino, J. B. Seoane--Sep\'{u}lveda, The
	Bohr radius of the $n$-dimensional polydisc is equivalent to $\sqrt{(\log
		n)/n}.$ \emph{Adv. Math.} \textbf{264} (2014), 726--746.
	
	\bibitem {boas}H. P. Boas, Majorant series. \emph{J. Korean Math. Soc.}
	\textbf{37} (2000), 321--337.
	
	\bibitem {bh}H. F. Bohnenblust, E. Hille, On the absolute convergence of
	Dirichlet series. \emph{Ann. of Math.} \textbf{32} (1931), 600--622.
	
	\bibitem {Defant}D. Carando, A. Defant, P. Sevilla--Peris, The
	Bohnenblust--Hille inequality combined with an inequality of Helson.
	\emph{Proc. Amer. Math. Soc.} \textbf{143} (2015), no. 12, 5233--5238.
	
	\bibitem{Cavalcante} W. Cavalcante, D. Pellegrino, Bohnenblust--Hille inequalities: analytical and computational aspects. 
	\emph{An. Acad. Brasil. Ciênc.} \textbf{91} (2019), no. 1, suppl. 1,  e20170398 (19 pages). 
	
	\bibitem {LLL}J. E. Littlewood, On bounded bilinear forms in an infinite
	number of variables. \emph{Quart. J.} (Oxford Ser.) \textbf{1} (1930), 164--174.
	
	\bibitem {Mariana}M. Maia, T. Nogueira, D. Pellegrino, The Bohnenblust--Hille
	inequality for polynomials whose monomials have a uniformly bounded number of
	variables. \emph{Integr. Equ. Oper. Theory} \textbf{88} (2017), 143--149.
	
	\bibitem {montanaro}A. Montanaro, Some applications of hypercontractive
	inequalities in quantum information theory. \emph{J. Math. Phys.} \textbf{53}
	(2012), no. 12, 122206 (15 pages).
	
	\bibitem{Nunez}D. N\'{u}\~{n}ez--Alarc\'{o}n, A note on the polynomial Bohnenblust--Hille inequality. \emph{J. Math. Anal. Appl.} \textbf{407}
	(2013), no. 1, 179--181.
	
	\bibitem {Teixeira}D. Pellegrino, E. Teixeira, Towards sharp Bohnenblust--Hille
	constantes. \emph{Commun. Contemp. Math.} \textbf{20} (2018), no. 3, 1750029
	(33 pages).
	
	\bibitem {Tomaz}D. Tomaz, Hardy--Littlewood inequalities for multipolynomials.
	\emph{Adv. Oper. Theory} \textbf{4} (2019), no. 3, 688--697.
	
	\bibitem {fv}F. Vieira Costa Junior, The optimal multilinear
	Bohnenblust--Hille constants: a computational solution for the real case.
	\emph{Numer. Func. Anal. Opt.} \textbf{39} (2018), 1656--1668.
	                                                        %
\end{thebibliography}
\end{document}